\documentclass[reqno]{proc-l}

\input supp-pdf.tex
\usepackage{latexsym}
\usepackage{amstext}
\usepackage {amsmath}
\usepackage {amsfonts}
\usepackage {amssymb}
\usepackage {amsthm}
\usepackage {bbm}
\usepackage{enumerate}
\usepackage{array}
\usepackage{comment}
\usepackage{xspace}
\newcommand{\eps}{\varepsilon}
\newcommand{\T}{T^\varepsilon_f}

\newcommand{\R}{\ensuremath{\mathbb{R}}}
\newcommand{\Rn}{\ensuremath{{\mathbb{R}^n}}}
\newcommand{\N}{\ensuremath{\mathbb{N}}\ }

\newcommand{\sphere}{{S^{n-1}}}

\newcommand{\al}{\alpha}
\newcommand{\ka}{\kappa}
\newcommand{\ep}{\varepsilon}
\def\Modd#1{\left\|#1\right\|}

\def\conv#1{#1\mathord{*}}

\newcounter{counter-liste}
\newcounter{counter-liste2}

  {\ \begin{list}{{(\arabic{counter-liste})}\hfill}%
  {\topsep2mm\itemindent1ex\leftmargin0cm\usecounter{counter-liste}}
  }%
  {\end{list}}
\newenvironment{liste(a)}%
  {\ \begin{list}{{(\alph{counter-liste})}\hfill}%
  {\topsep2mm\itemindent1ex\leftmargin0cm\usecounter{counter-liste}}
  }%
  {\end{list}}
\def\R{\mathbb R}

\let\e\varepsilon

\def\pref#1{(\ref{#1})}
\theoremstyle{plain}
\numberwithin{equation}{section}
\newtheorem{lemma}{Lemma}[section]
\newtheorem{theorem}[lemma]{Theorem}

\newtheorem{corollary}[lemma]{Corollary}
\newtheorem{problem}[lemma]{Problem}
\theoremstyle{definition}
\newtheorem{remark}[lemma]{Remark}
\begin{document}
\title{Sobolev regularity and an enhanced Jensen inequality}

\author[Peletier]{Mark A. Peletier}
\address{Technische Universiteit Eindhoven,
Den Dolech 2,
P.O. Box 513,
5600 MB Eindhoven}
\thanks{MA and MR were supported by NWO grant 639.032.306. 
RP wishes to thank the Centrum voor Wiskunde en Informatica for their 
financial support.}

\author[Planqu\'{e}]{Robert Planqu\'{e}}
\address{Department of Mathematics, 
Vrije Universiteit,
De Boelelaan 1081a,
1081 HV Amsterdam}
\email{rplanque@few.vu.nl}

\author[R\"{o}ger]{Matthias R\"{o}ger}
\address{Max Planck Institute for Mathematics in the Sciences,
Inselstr. 22,
D-04103 Leipzig}



\date{\today}

\begin{abstract}
We derive a new criterion for a real-valued 
function $u$ to be in the Sobolev space
$W^{1,2}(\R^n)$. This criterion consists of comparing the value of a
functional $\int f(u)$ with the values of the same functional applied to 
convolutions of $u$ with a Dirac sequence. The difference of these
values converges to zero as the convolutions approach $u$, and we
prove that the rate of convergence to zero is connected to
regularity: $u\in W^{1,2}$ if and only if the convergence is sufficiently fast. We
finally apply our criterium to a minimization problem with constraints,
where regularity of minimizers cannot be deduced from the 
Euler-Lagrange equation.
\end{abstract}

\maketitle

\section{Introduction}
Jensen's inequality states that if $f:\R\to\R$ is convex and
$\varphi\in L^1(\R^n)$ with $\varphi\geq 0$ and $\int\varphi = 1$,
then
\[
\int_{\R^n} f\bigl(u(x)\bigr) \varphi(x)\, dx
\geq 
f\left(\int_{\R^n} u(x)\varphi(x)\, dx\right) ,
\]
for any $u:\R^n\to\R$ for which the integrals make sense. A consequence of this inequality is that 
\begin{gather}
  \int_\Rn f(u)\,=\, \int_\Rn f(u)*\varphi\,\geq\, \int_\Rn
  f(u*\varphi).
 \label{eq:motiv-jensen}
\end{gather}

In this paper we investigate the inequality~\pref{eq:motiv-jensen}
more closely. In particular we study the relationship between the
regularity of a function $u\in L^2(\R^n)$ and the asymptotic behaviour
as $\eps\to0$ of
\begin{gather}
  \T(u)\,:=\,\int_\Rn \bigl[f(u) - f(u*\varphi_\eps)\bigr].\label{eq:ex-main}
\end{gather}
Here $\varphi_\eps(y)\,=\, \eps^{-n}\varphi(\eps^{-1}y)$ and $f$ is a
smooth function, but now not necessarily convex. Since
$u*\varphi_\eps\to u$ almost everywhere, we find for `well-behaved'
$f$ that $\T(u)\to0$ as $\eps\to0$. Our aim is to establish a
connection between the \emph{rate of convergence} of $\T(u)$ to zero
and the \emph{regularity} of~$u$.

Such a connection between the decay rate of $\T(u)$ and the regularity
of $u$ is suggested by the following informal arguments. Taking $n=1$,
and assuming $\varphi$ to be even and $u$ to be smooth, we develop
$u*\varphi_\eps$ as
\begin{gather}
  u*\varphi_\eps(x)\,=\, \int_\R u(x-y)\varphi_\eps(y)\,dy\,\approx\,
  u(x) +
  \frac{\eps^2}{2}u^{\prime\prime}(x) \int_\R  y^2\varphi(y)\,dy,
\end{gather}
so that
\begin{gather}
 \T(u)   \,\approx\, -c\eps^2\int_\R f^\prime(u)
 u^{\prime\prime}\,=\, c\eps^2 \int_\R
 f^{\prime\prime}(u)|u^\prime|^2, \label{eq:formal} 
\end{gather}
with $c=\frac{1}{2}\int_\R  y^2\varphi(y)\,dy$.
This suggests that if $u\in W^{1,2}(\R)$ then $\T(u)$
is of order~$\eps^2$. Moreover, since
$f^{\prime\prime}\geq 0$ in the case that $f$ is convex,
\eqref{eq:formal} gives an enhanced version of
\eqref{eq:motiv-jensen}. 

Conversely, for functions not in $W^{1,2}(\Rn)$ we might not observe
decay of $\T(u)$ with the rate of $\eps^2$, as a simple example
shows. Take $f(u) = u^2$ and consider a function with a jump
singularity such as $u(x) = H(x) - H(x-1)\not\in W^{1,2}(\R)$, where
$H(x)$ is the Heaviside function. Now choose as regularization kernels
the functions $\varphi_\eps (x) = \frac 1{2\eps}(H(x+\eps) -
H(x-\eps))$, after which an explicit calculation shows that $\T(u) =
2\eps/3$. Here the decay is only of order $\eps$.

The goal of this paper is to prove the asymptotic development
\eqref{eq:formal} in arbitrary dimension and to show that also the
converse statement holds true: if $\T(u)$ is of order $\eps^2$ then
$u\in W^{1,2}(\Rn)$. These results are stated below and proved in
Sections~\ref{sec:proof1} and~\ref{sec:proof2}.

The fact that one can deduce regularity from the decay rate is the
original motivation of this work. In Section~\ref{sec:application} we
illustrate the use of this result with a minimization problem in which
the Euler-Lagrange equation provides no regularity for a minimizer
$u$. Instead we estimate the decay of $\T(u)$ directly from the
minimization property and obtain that $u$ is in $W^{1,2}$.

\medskip
Our results can be compared to the characterisation of Sobolev spaces
introduced by Bourgain, Brezis and Mironescu \cite{BBM,Pon}. In fact,
the regularity conclusion in Theorem \ref{the:reg} could be derived
from \cite{BBM}, with about the same amount of effort as the
self-contained proof that we give here.

\subsection{Notation and assumptions}\label{sec:ass}
Let $\varphi\in L^1(\R^n)$ satisfy
\begin{align}
  &\varphi\,\geq\, 0\quad\text{ in }\R^n,  \label{ass:phi-pos}\\
  & \int_\Rn \varphi\,=\, 1,\label{ass:phi-mass}\\
  &\!\!\int_\Rn y\varphi(y)\,dy\,=\, 0,\label{ass:phi-asym}\\
  &\!\!\int_\Rn |y|^2\varphi(y)\,dy\,<\,\infty. \label{ass:2-moment}
\end{align}
For $\eps>0$ we define the Dirac sequence
\begin{gather}\label{eq:def-phi-eps}
  \varphi_\eps(x)\,:=\, \frac{1}{\eps^{n}}\varphi\Big(\frac{x}{\eps}\Big).
\end{gather}
Eventually we will restrict ourselves to the case that $\varphi$ is
rotationally symmetric, that is $\varphi(x)=\tilde{\varphi}(|x|)$, where
$\tilde{\varphi}:\R^+_0\to\R^+_0$. Then \eqref{ass:phi-mass} is equivalent to
\begin{gather}
  n\omega_n\int_0^\infty
  r^{n-1}\tilde{\varphi_\eps}(r)\,dr\,=\, 1, \label{eq:kernel-1}
\end{gather}
where $\omega_n$ denotes the volume of the unit-ball, i.e.
$|B_1(0)|=\omega_n$.

For a function $u\in L^1(\R^n)$ and $0\leq s\leq 1$ we define the
convolution $u_\eps$ and modified convolutions $u_{\eps,s}$ by
\begin{align}
  u_\eps(x)\,&:=\, \big(u*\varphi_\eps\big)(x)\,=\, \int_\Rn
  u(x-y)\varphi_\eps(y)\,dy,\label{eq:def-conv}\\
  u_{\eps,s}(x)\,&:=\, \big(u*\varphi_{\eps s}\big)(x)
      \,=\, \int_\Rn u(x-y)\varphi_{\eps s}(y)\,dy
    \,=\, \int_\Rn u(x-sz)\varphi_\eps(z)\,dz.\label{eq:def-mod-conv}
\end{align}
Note that $u_\eps\,=\, u_{\eps,1}, u\,=\, u_{\eps,0}.$

We also use the notation $a\otimes b$ for the tensor product of $a,b\in\Rn$.

\subsection{Statement of main results}
Our first result proves \eqref{eq:formal}.
\begin{theorem}\label{the:main}
Let $f\in C^2(\R)$ have uniformly bounded second derivative and
\begin{gather}\label{eq:ass-f}
  f(0)=0,\qquad f^\prime(0)\,=\,0.
\end{gather}
Let $(\varphi_\eps)_{\eps>0}$ be a
Dirac sequence as in Section \ref{sec:ass}.
Then for any $u\in W^{1,2}(\Rn)$,
\begin{gather}\label{eq:the-main}
  \lim_{\eps\to 0}\frac{1}{\eps^2}\int_\Rn \bigl[f(u)-f(u_\eps)\bigr]\,dx
  =\, \frac{1}{2}\int_\Rn f^{\prime\prime}(u(x))\nabla u(x)\cdot A(\varphi)
  \nabla u(x)\,dx,
\end{gather}
where 
\begin{gather}
  A(\varphi)\,=\, \int_\Rn \big(y\otimes
  y\big)\varphi(y)\,dy.
\end{gather}
If $\varphi$ is rotationally symmetric, i.e.
$\varphi(x)=\tilde{\varphi}(|x|)$, then
\begin{gather}\label{eq:the-main-2}
  \lim_{\eps\to 0}\frac{1}{\eps^2}\int_\Rn \bigl[f(u)-f(u_\eps)\bigr]\,dx
  =\, \frac{1}{2} a(\varphi)\int_\Rn f^{\prime\prime}(u(x))|\nabla u(x)|^2\,dx,
\end{gather}
with
\begin{gather}
  a(\varphi)\,=\, \omega_n\int_0^\infty
  r^{n+1} \tilde{\varphi}(r)\,dr.
\end{gather}
\end{theorem}
The second theorem shows that for uniformly convex $f$ a decay of $\T(u)$
of order~$\eps^2$ implies that $u \in W^{1,2}(\Rn)$.
\begin{theorem}\label{the:reg}
Let $f\in C^2(\R)$ have uniformly bounded second derivative, 
assume that \eqref{eq:ass-f} is satisfied, and that there is a positive
number $c_1>0$ such that 
\begin{gather}\label{ass:strictly-convex}
  f''\geq c_1\qquad \text{on $\R$.}
\end{gather}
If for $u\in L^2(\Rn)$,
\begin{gather}
  \liminf_{\eps\to 0} \frac 1{\eps^2}\int_\Rn \bigl[f(u(x))-f(u_\eps(x))\bigr]\,dx\,<\,\infty,
\end{gather}
then $u\in W^{1,2}(\Rn)$. In particular \eqref{eq:the-main} and
\eqref{eq:the-main-2} hold.
\end{theorem}
\begin{remark}\label{rem:ass}
We prescribe \eqref{eq:ass-f} since in general the difference
$f(u)-f(u_\eps)$ need not have sufficient decay to 
be Lebesgue integrable. For functions $u\in L^1(\Rn)\cap L^2(\Rn)$ Theorems~\ref{the:main} and~\ref{the:reg} hold
even without assuming \eqref{eq:ass-f}.
\end{remark}

\section{Proof of Theorem \ref{the:main}}
\label{sec:proof1}
We first show that $f(u)-f(u_\eps)\in L^1(\Rn)$. Let $\eta\in
C^0_c(\Rn;[0,1])$. 
Using the Fundamental Theorem of Calculus we obtain that
\begin{align}
  \lefteqn{\int_\Rn \eta(x)\bigl| f(u(x))-f(u_\eps(x))\bigr|\,dx = }
       \qquad\qquad&\notag\\
  &=\, \int_\Rn \eta(x)\Big|\int_0^1\frac{\partial}{\partial s}
  f\big(u_{\eps,s}(x)\big) \,ds\Big|\,dx\notag\\
  &\leq\, \int_0^1\int_\Rn \eta(x)\int_\Rn \Big| f^\prime\big(u_{\eps,s}(x)\big)\nabla
  u(x-sy)\cdot y\varphi_\eps(y)\Big|\,dy\,dx\,ds\notag\\
  &\leq\, \frac{1}{2}\int_0^1\int_\Rn \eta(x)\int_\Rn f^\prime\big(u_{\eps,s}(x)\big)^2
  |y|\varphi_\eps(y)\,dy\,dx\,ds\notag\\
  &\qquad +\frac{1}{2}\int_0^1\int_\Rn \eta(x) \int_\Rn |\nabla
  u(x-sy)|^2 |y|\varphi_\eps(y)\,dy\,dx\,ds\label{eq:L1-1},
\end{align}
where we have used Fubini's Theorem and Young's inequality.
For the first term on the right-hand side we deduce by
\eqref{ass:2-moment} and \eqref{eq:ass-f} that
\begin{align}
  \int_\Rn \int_\Rn\eta(x)
       f^\prime\big(u_{\eps,s}(x)\big)^2 |y|\varphi_\eps(y)\,dy\,dx 
  &\;\leq\; C_\eps(\varphi)\|f''\|_\infty^2\int_\Rn u_{\eps,s}(x)^2\,dx\notag\\
  &\stackrel{\pref{eq:motiv-jensen}}\leq C_\eps(\varphi)\|f''\|_\infty^2\int_\Rn u(x)^2\,dx. \label{eq:L1-2}
\end{align}
For the second term on the right-hand side of \eqref{eq:L1-1} we obtain
similarly 
\begin{align}
  \int_\Rn \int_\Rn\eta(x) |\nabla
  u(x-sy)|^2 |y|\varphi_\eps(y)\,dy\,dx\,\leq\, C_\eps(\varphi)\int
  |\nabla u(x)|^2\,dx.\label{eq:L1-3}
\end{align}
By \eqref{eq:L1-1}-\eqref{eq:L1-3} we therefore obtain
\begin{align}
  \int_\Rn \eta(x)\bigl| f(u(x))-f(u_\eps(x))\bigr|\,dx
   \,\leq\,
   C_\eps(\varphi)\big(1+\|f''\|_\infty^2\big)\|u\|_{W^{1,2}(\Rn)}^2\,<\,\infty. 
\end{align}
Letting $\eta\nearrow 1$ we deduce that $f(u)-f(u_\eps)\in L^1(\Rn)$.
Repeating some of the
calculations above, we obtain that
\begin{align}
  \T(u)\,&=\, -\int_\Rn \int_0^1\frac{\partial}{\partial s}
  f\big(u_{\eps,s}(x)\big) \,ds\,dx\notag\\
  &=\, \int_0^1\int_\Rn \int_\Rn  f^\prime\big(u_{\eps,s}(x)\big)\nabla
  u(x-sy)\cdot y\varphi_\eps(y)\,dy\,dx\,ds \label{eq:proof1}
\end{align}
and that
\begin{align}
  f^\prime\big(u_{\eps,s}(x)\big) - f^\prime\big(u_{\eps,s}(x-sy)\big) 
  \,&=\, -\int_0^s
  \frac{\partial}{\partial r} f^\prime\big(u_{\eps,s}(x-ry)\big)\,dr
  \notag\\
  &=\, \int_0^s f^{\prime\prime}\big(u_{\eps,s}(x-ry)\big)\nabla u_{\eps,s}(x-ry)\cdot y\,dr.
  \label{eq:proof2}
\end{align}
Since for all $0\leq s\leq 1$
\begin{align*}
  \lefteqn{\int_\Rn\int_\Rn  f^\prime\big(u_{\eps,s}(x-sy)\big)\nabla
  u(x-sy)\cdot y\varphi_\eps(y)\,dy\,dx \,=\,}  
\qquad\qquad\qquad&\\
  &\;\,=\,\int_\Rn\int_\Rn  f^\prime\big(u_{\eps,s}(z)\big)\nabla
  u(z)\cdot y\varphi_\eps(y)\,dz\,dy\\
  &\stackrel{\pref{ass:phi-asym}}=\;0,
\end{align*}
we deduce from \eqref{eq:proof1}, \eqref{eq:proof2} that
\begin{align}
  \lefteqn{\frac{1}{\eps^2}\int_\Rn \bigl[f(u(x))-f(u_\eps(x))\bigr]\,dx =}\qquad\qquad\notag\\
  =\,& \frac{1}{\eps^2}\int_0^1\int_\Rn \int_\Rn \Big[\Big(\int_0^s
  f^{\prime\prime}\big(u_{\eps,s}(x-ry)\big)
  \nabla u_{\eps,s}(x-ry)
  \cdot y\,dr\Big)  \notag\\
  &\qquad\qquad\qquad\qquad
  \nabla u(x-sy)\cdot y\varphi_\eps(y)\Big] dy\,dx\,ds\notag\\
  =\,& \frac{1}{\eps^2}\int_0^1\int_0^s\int_\Rn \int_\Rn \Big[
  f^{\prime\prime}\big(u_{\eps,s}(x)\big)
  \nabla u_{\eps,s}(x)
  \cdot y\notag\\
  &\qquad\qquad\qquad\qquad
  \nabla u(x-ry)\cdot y\varphi_\eps(y)\Big] dy\,dx\,dr\,ds\notag\\
  =\,& \frac{1}{\eps^{n+2}}\int_0^1\int_0^s\int_\Rn
  f^{\prime\prime}\big(u_{\eps,s}(x)\big)
  \nabla u_{\eps,s}(x)
  \cdot{} \notag\\
  &\qquad\qquad\qquad\qquad
  \int_\Rn \bigl(y \otimes y\bigr) \cdot\nabla u(x-ry)\,\varphi(\eps^{-1}y)\, dy\,dx\,dr\,ds\notag\\
  =\,& \frac{1}{\eps^{n+2}}\int_0^1\int_0^s\int_\Rn
  f^{\prime\prime}\big(u_{\eps,s}(x)\big)
  \nabla u_{\eps,s}(x)
  \cdot{} \notag\\
  &\qquad\qquad\qquad\qquad
  \frac1{r^{n+2}}\int_\Rn \bigl(z \otimes z\bigr) \cdot\nabla
  u(x-z)\,\varphi(\eps^{-1}r^{-1}z)\, dz\,dx\,dr\,ds\notag\\ 
  =\,& \int_0^1\int_0^s \int_\Rn f^{\prime\prime}\big(u_{\eps,s}(x)\big) \nabla
  u_{\eps,s}(x) \cdot \big(\kappa_{\eps r}* \nabla u\big)(x)\,dx\,dr\,ds,
 \label{eq:proof3} 
\end{align}
with the modified convolution kernel $\kappa_{\eps r}$ defined by
\begin{gather}
  \kappa_{\eps r}(z)\,:=\, (\eps r)^{-n}\kappa\Big(\frac{z}{\eps r}\Big),
\end{gather}
where $\kappa$ is given by
\begin{gather}
  \kappa(y)\,:=\, \big(y\otimes
  y\big)\varphi(y).
\end{gather}
As $\eps\to 0$ we have that for all $0\leq r$, $s\leq 1$,
\begin{alignat}{2}
  u_{\eps,s}\,&\to\, u\quad&&\text{ almost everywhere,}\notag\\
  \nabla u_{\eps,s}\,&\to\, \nabla u\quad&&\text{ in }L^2(\R^n),
    \label{conv:nabla-u}\\
  \kappa_{\eps r}* \nabla u\,&\to\, \Big(\int_\Rn
  \kappa(y)\,dy\Big)\nabla u\quad&&\text{ in }L^2(\R^n).
    \label{conv:kappa-u}
\end{alignat}
Moreover, the integrand on the right-hand side of \eqref{eq:proof3} is
dominated by the function
\begin{gather}
  \frac{1}{2}\|f^{\prime\prime}\|_{L^\infty(\R)}\Big( |\nabla u_{\eps,s}|^2 +
  \big|\kappa_{\eps r}* \nabla u\big|^2 \Big),
\end{gather}
which converges strongly in $L^1(\R^n)$ as $\eps\to 0$
by~\pref{conv:nabla-u} and~\pref{conv:kappa-u}.  
We therefore deduce from \eqref{eq:proof3} that
\begin{align}
  \lefteqn{\lim_{\eps\to 0}\frac{1}{\eps^2}\int_\Rn
  \bigl[f(u(x))-f(u_\eps(x))\bigr]\,dx = }\qquad\qquad\notag\\ 
  =\,& \int_0^1\int_0^s \int_\Rn f^{\prime\prime}(u(x)) \nabla
  u(x)\cdot\Big(\int_\Rn
  \kappa(y)\,dy\Big)\nabla u(x)\,dx\,dr\,ds\notag\\
  =\,& \frac{1}{2}\int_\Rn f^{\prime\prime}(u(x))\nabla u(x)\cdot\Big(\int_\Rn
  \kappa(y)\,dy\Big)\nabla u(x)\,dx.
\end{align}
In the rotationally symmetric case we observe that
\begin{align}
  \int_\Rn  \big(y\otimes  y \big)\varphi(y)\,dy\,&=\, \Big(\int_0^\infty
  r^{n+1} \tilde{\varphi}(r)\,dr\Big)\int_\sphere
  \theta\otimes\theta\,d\theta\notag\\
  &=\, \omega_n\Big(\int_0^\infty
  r^{n+1} \tilde{\varphi}(r)\,dr\Big)Id.\label{eq:proof5}
\end{align}
Here we used that for $i,j=1,\ldots,n$,
\begin{gather*}
  \int_\sphere \theta_i \theta_j \,d\theta\,=
  \begin{cases}
    \omega_n &\text{ if }i=j,\\
    0 &\text{ if }i\neq j.
  \end{cases}
\end{gather*}  
We therefore obtain from \eqref{eq:proof3}-\eqref{eq:proof5} that in the
rotationally symmetric case
\begin{align}
  \lim_{\eps\to 0}\frac{1}{\eps^2}\int_\Rn \bigl[f(u)-f(u_\eps)\bigr]\,dx\,=\, 
  \omega_n\Big(\int_0^\infty
  r^{n+1} \tilde{\varphi}(r)\,dr\Big)\frac{1}{2}\int_\Rn f^{\prime\prime}(u)
  |\nabla u|^2\,dx.
\end{align}

\section{Proof of Theorem \ref{the:reg}}
\label{sec:proof2}
Let
\begin{gather}
  \Lambda\,:=\, \liminf_{\eps\to 0} \frac1{\eps^2}\int_\Rn \bigl[f(u(x))-f(u_\eps(x))\bigr]\,dx.
\end{gather}
We first remark that by assumption \eqref{ass:strictly-convex} the
function $r\mapsto f(r)-\frac{c_1}{2}r^2$ is convex and we deduce again by~\pref{eq:motiv-jensen} that
\begin{gather}
  \int_\Rn \Bigl[\Big(f(u) -\frac{c_1}{2} u^2\Big) - \Big(f(u_\eps) - \frac{c_1}{2}
  u_\eps^2\Big)\Bigr]\,\geq\, 0.
\end{gather}
This implies that
\begin{gather}
  \liminf_{\eps\to 0}\frac{1}{\eps^2} \int_\Rn \Big(u^2 - u_\eps^2\Big)
  \,\leq\, 
  \frac{2}{c_1}\liminf_{\eps\to 0} \frac{1}{\eps^2}\int_\Rn \bigl[ f(u)
  -f(u_\eps)\bigr]\,=\, \frac{2}{c_1}\Lambda.\label{eq:est-squares} 
\end{gather}
Next we consider $\delta>0$ and the regularizations $u_\delta=
u*\varphi_\delta$. Set $\gamma_\eps:=\varphi_\eps*\varphi_\eps$.
Then we obtain that
\begin{align}
  \lefteqn{\int_\Rn \bigl[u_\delta^2(x)
  -\big(u_\delta*\varphi_\eps\big)^2(x)\bigr]\,dx =
  }\qquad\qquad\notag\\ 
  =\;\;& \int_\Rn u_\delta^2 -\int_\Rn\int_\Rn u_\delta(x)u_\delta(y)
  \gamma_\eps(x-y) \,dy\,dx\notag\\
  =\;\;& \frac{1}{2}\int_\Rn\int_\Rn \big( u_\delta(x)-
  u_\delta(y)\big)^2\gamma_\eps(x-y)\,dy\,dx \notag\\
  \stackrel{\pref{eq:motiv-jensen}}\leq\,& \frac{1}{2}\int_\Rn\int_\Rn\int_\Rn
  \big(u(x-z)-u(y-z)\big)^2\varphi_\delta(z)\, dz\,\gamma_\eps(x-y)\,dy\,dx \notag\\
  =\;\;& \frac{1}{2}\int_\Rn\int_\Rn\int_\Rn
  \big(u(\xi)-u(\eta)\big)^2\gamma_\eps(\xi-\eta)\,d\eta\,d\xi\,\varphi_\delta(z)\,dz
  \notag\\ 
  =\;\;& \frac{1}{2}\int_\Rn\int_\Rn \big( u(\xi)-
  u(\eta)\big)^2\gamma_\eps(\xi-\eta)\,d\eta\,d\xi \notag\\
  =\;\;&\int_\Rn \big[u^2(\xi) -u_\eps^2(\xi)\bigr]\,d\xi.\label{eq:est-diff-squares}
\end{align}
We therefore deduce from Theorem \ref{the:main} and
\eqref{eq:est-squares}, \eqref{eq:est-diff-squares} that
\begin{align}
  \int_\Rn |\nabla u_\delta|^2\,&=\, C(\varphi)\liminf_{\eps\to
  0}\frac1{\eps^2}\int_\Rn \bigl[u_\delta^2(x)
  -\big(u_\delta*\varphi_\eps\big)^2(x)\bigr]\,dx\notag\\ 
  &\leq\, C(\varphi)\liminf_{\eps\to
  0}\frac1{\eps^2}\int_\Rn \bigl[u^2(x) -u_\eps^2(x)\bigr]\,dx\,\leq\, \frac{2}{c_1}\Lambda
  C(\varphi). 
\end{align}
Since this estimate is uniform in $\delta>0$ it follows that $u\in
W^{1,2}(\Rn)$. By Theorem~\ref{the:main} we therefore deduce~\eqref{eq:the-main} and~\eqref{eq:the-main-2}.

\section{Application: regularity of minimizers in a lipid bilayer model} 
\label{sec:application}

We illustrate the utility of Theorems~\ref{the:main} and~\ref{the:reg}
with an example. The problem is to determine the regularity of solutions
of the following minimization problem.
\begin{problem}
Let $\alpha>0, h>0$ and a  kernel $\kappa\in W^{1,1}(\R,\R^+_0)$ be given
such that $\kappa^\prime\in BV(\R)$. Denote by $\tau_h$ the
translation operator, 
\[
(\tau_h u)(x) := u(x-h)\quad\text{ for }u:\R\to\R.
\]
Consider the set $K\subset\R$,
\begin{gather}
   K \,:=\, \Big\{u \in L^1(\R)\ \Big|\ \int_\R u =1,\ u \ge 0,\ u + \tau_h u
  \le 1 \text{ a.e.}\Big\}\label{eq:def-K}
\end{gather}
and a strictly convex, increasing, and smooth function $f$ with
$f(0)=f'(0)=0$. We then define a functional $F:K\to\R$ by
\begin{align}
  F(u) &:= \int_\R f(u) - \int_\R \al u \conv \ka u \label{eq:def-F}
\end{align}
and consider the problem of finding a function $u \in K$ that minimizes $F$ in $K$, \emph{i.e.}
\begin{equation}
\label{MIeq:theproblem}
F(u) = \min\{F(v)\ |\ v \in K\}.
\end{equation}
\end{problem}

\begin{remark}
This problem arises in the modelling of lipid bilayers, biological
membranes (see~\cite{BlomPeletier04} for more details).
There, a specific type of molecules is considered, consisting of two \emph{beads}
connected by a rigid rod. The beads have a certain volume, but the rod
occupies no space. 

The beads are assumed to be of sub-continuum size, but the rod length is
non-negligeable at the continuum scale. In the one-dimensional case we
also assume that the rods lie parallel to the single spatial
axis. Combining these assumptions we model the distribution of such
molecules over the real line by a variable $u$ that represents the
volume fraction occupied by \emph{leftmost} beads. The volume fraction
of rightmost beads is then given by $\tau_h u$, and in the condition
$u+\tau_h u\leq1$ we now recognize a volume constraint. 
 
The functional $F$ represents a free energy. In the entropy term $\int
f(u)$ the function $f$ is strictly convex, increasing,  and smooth and
can be assumed to satisfy $f(0)=f'(0)=0$. The destabilizing term $-\int
u\conv \ka u$ is a highly stylized representation of the hydrophobic
effect, which favours clustering of beads. 
\end{remark}

Without the constraints $u\geq 0$ and $u + \tau_h u \le 1$ present in
\eqref{eq:def-K}, we 
would immediately be able to infer that minimizers (even stationary
points) are
smooth using a simple bootstrap argument: since $u\in K$ we have $u\in
L^p(\R)$ and $\conv \ka u \in W^{1,p}$ for all $1\leq p\leq\infty$, and
therefore, using the Euler-Lagrange 
equation, 
\[
f'(u) - 2\al \conv \ka u = \lambda,
\]
for some $\lambda \in \R$, we find $u \in W^{1,p}(\R)$. Iterating this
procedure we obtain that $u\in  W^{k,p}(\R)$ for all $k\in\N$, $1\leq p
\leq\infty$. 
However, when including the two constraints, two additional Lagrange
multipliers $\mu$ and $\nu$ appear in
the Euler-Lagrange equation,  
\[
f'(u) - 2\al \conv \ka u = \lambda +\mu - \nu - \tau_{-h}\nu.
\]
Here $\mu$ and $\nu$ are measures on $\R$, and standard theory provides no further regularity than this. In this case the lack of regularity in the
right-hand side interferes with the bootstrap process, and this equation
therefore does not give rise to any additional regularity. 

\medskip
The interest of Theorem~\ref{the:reg} for this case lies in the fact
that $K$ is closed under convolution. In particular, if $u$ is a minimizer of $F$ in
$K$, then $u_\e = \conv{\varphi_\eps} u$ is also admissible, and we can
compare $F(u_\e)$ with $F(u)$. From this comparison and an
application of Theorem \ref{the:reg} we deduce the
regularity of~$u$:

\begin{corollary}
\label{cor:membranes} Let $u$ be a minimizer
of~\eqref{MIeq:theproblem}. Then $u \in W^{1,2}(\R)$. 
\end{corollary}

\begin{proof}
Choose a function $\varphi \in L^1(\R)$ as in Section~\ref{sec:ass} with
Dirac sequence 
$(\varphi_\ep)_{\ep>0}$, and set $u_\ep := \conv
{\varphi_\ep}u$. First we show that there exists a constant $C\in \R$
such that for all $u\in L^2(\R)$
\begin{equation}
  \label{MEMeq:rho_ep*e_order_ep^2}
  \Big|\int_\R\big( u\conv \ka u -
  u_\ep\conv\ka u_\ep\big)\Big|\,\le\, C \eps^2 \int_\R u^2.
\end{equation}
With this aim we observe that 
\begin{gather}
  \int_\R\big( u\conv \ka u -
  u_\ep\conv\ka u_\ep\big)\,=\, \int_\R u\big(\conv\kappa u - \conv{\conv
  \kappa u}{\gamma_\eps}\big) \label{eq:gamma-eps}
\end{gather}
with
\begin{align*}
  \gamma\,&:=\,\conv{\varphi}{\varphi},\qquad
  \gamma_\eps(x)\,:=\,\big(\conv{\varphi_\eps}{\varphi_\eps}\big)(x)
  \,=\,\frac{1}{\eps}\gamma\Big(\frac{x}{\eps}\Big),
\end{align*}
and that $\gamma,\gamma_\eps$ satisfy the assumptions in Section \ref{sec:ass}.

Since $\kappa'\in BV(\R)$ we have $(\kappa*u)''=\kappa''*u\in L^2(\R)$
and
\begin{gather*}
  \|\kappa''*u\|_{L^2(\R)}\,\leq\, \Big(\int_\R
    |\kappa''|\Big)\|u\|_{L^2(\R)}. 
\end{gather*}
Repeating some arguments of Theorem \ref{the:main} we calculate
that 
\begin{align*}
  (\conv\kappa u)(x) - (\conv\kappa\conv u\gamma_\eps)(x)&=\,
    - \int_0^1 \frac d{ds} \int_\R (\conv\kappa u)(x-sy)\gamma_\eps(y)
       \,dy\,ds\,\\
  &=\, \int_0^1\int_\R (\conv\kappa u)'(x-sy)y\gamma_\eps(y)\,dy\,ds\\
  &=\, \int_0^1 \int_0^1 \int_\R (\conv\kappa u)''(x-rsy) sy^2 \gamma_\eps(y)
     \,dy\,ds\,dr \\
  &=\, \eps^2 \int_0^1\int_0^1\int_\R (\conv\kappa u)^{\prime\prime}(x-\eps
  rsz)sz^2\gamma(z)\,dz\,dr\,ds ,
\end{align*}
and therefore that
\begin{align*}
  &{\Big|\int_\R u\big(\conv\kappa u - \conv{\conv
  \kappa u}{\gamma_\eps}\big)\Big|\,\leq\,}\quad\quad\\
  \leq\, &\eps^2 \int_0^1\int_0^1\int_\R \Big[ \int_\R (\conv\kappa u)^{\prime\prime}(x-\eps
  rsz)^2\,dx\Big]^{\frac{1}{2}}\Big[\int_\R
  u(x)^2\,dx\Big]^{\frac{1}{2}}sz^2\gamma(z)\,dz\,dr\,ds\\
  \leq\, &\eps^2 C(\gamma)\|\kappa''*u\|_{L^2(\R)}\|u\|_{L^2(\R)}\\
  \leq\, &\eps^2 C(\varphi,\kappa)\|u\|_{L^2(\R)}^2,
\end{align*}
which by \eqref{eq:gamma-eps} implies \eqref{MEMeq:rho_ep*e_order_ep^2}.
%

Since $K$ is closed under convolution with $\varphi_\ep$, $u_\ep$ is
admissible, and therefore
\[
0\leq F(u_\ep)-F(u) \leq \int_\R [f(u_\ep) - f(u)] + C \ep^2 \int_\R u^2.
\]
The last term is bounded by
\[
C\ep^2\Modd u_{L^1(\R)}\Modd u_{L^\infty(\R)}\leq C \ep^2,
\]
which follows from combining both inequality constraints in
\eqref{eq:def-K}.  Therefore $\int [f(u)-f(u_\ep)] = O(\ep^2)$, and
from Theorem~\ref{the:reg} we conclude that $u \in W^{1,2}(\R)$.

\end{proof}

As is described in more detail in the thesis~\cite{Planque05}, Corollory~\ref{cor:membranes} paves the way
towards rigorously deriving the Euler-Lagrange equations for
Problem~\pref{MIeq:theproblem}, and forms an 
important ingredient of the  proof of existence of minimizers.


\providecommand{\bysame}{\leavevmode\hbox to3em{\hrulefill}\thinspace}
\providecommand{\MR}{\relax\ifhmode\unskip\space\fi MR }
\providecommand{\MRhref}[2]{%
  \href{http://www.ams.org/mathscinet-getitem?mr=#1}{#2}
}
\providecommand{\href}[2]{#2}

\end{document}